\documentclass[journal]{IEEEtran}
\renewcommand{\baselinestretch}{.94}
\usepackage{graphicx}
\usepackage{epstopdf}
\usepackage{dsfont}
\usepackage{amsmath}
\usepackage[hidelinks]{hyperref}
\usepackage{bm}
\usepackage{url}
\usepackage{caption}
\usepackage{subcaption}
\usepackage{wrapfig}
\usepackage{lastpage}
\usepackage{fancyhdr}
\usepackage{amssymb}
\usepackage{cite}
\usepackage{hyphenat}
\title{A Novel \textit{a priori} State Computation Strategy for the Unscented Kalman Filter to Improve Computational Efficiency}

\author{\IEEEauthorblockN{Sanat K. Biswas\IEEEauthorrefmark{1},
Li Qiao\IEEEauthorrefmark{2},
Andrew G. Dempster\IEEEauthorrefmark{1}}

\IEEEauthorblockA{\IEEEauthorrefmark{1}
Australian Centre for Space Engineering Research,
The University of New South Wales, NSW, Australia, 2052}

\IEEEauthorblockA{\IEEEauthorrefmark{2}
School of Engineering and Information Technology,
The University of New South Wales, Canberra, ACT, Australia, 2600}}

\begin{document}
\maketitle
\begin{abstract}
A priori state vector and error covariance computation for the Unscented Kalman Filter (UKF) is described. The original UKF propagates multiple sigma points to compute the \textit{a priori} mean state vector and the error covariance, resulting in a higher computational time compared to the Extended Kalman Filter (EKF). In the proposed method, the posterior mean state vector is propagated and then the sigma points at the current time step are calculated using the first-order Taylor Series approximation. This reduces the computation time significantly, as demonstrated using two example applications which show improvements of 90.5\% and 92.6\%. This method shows the estimated state vector and the error covariance are accurate to the first-order Taylor series terms. A second method using Richardson Extrapolation improves prediction accuracy to the second-order Taylor series terms. This is implemented on the two examples, improving efficiency by 85.5\% and 86.8\%.
\end{abstract}
\section{Introduction} 
The Kalman Filter is a long-established state estimation technique and many variants of this filter are widely used in numerous engineering applications, for example: object tracking, navigation, computer vision, economics and many more. The classical Kalman Filter was designed to address the estimation problem for linear systems \cite{Kalman1960}. NASA Ames Research Centre applied this optimal estimation formula for estimating the position and the velocity of a space vehicle. As the dynamics of a space vehicle are non-linear, the system was linearised using a first-order Taylor series approximation around an operating region to calculate the conditional error covariance and the Kalman gain \cite{Smith1964,Bernstein1966}. This is known as the Extended Kalman Filter (EKF) and is now widely used for numerous non-linear estimation problems. However this approach leads to a suboptimal solution to the non-linear estimation problem \cite{Cox1964,Athans1968} and requires an additional process noise covariance matrix for convergence of the solution \cite{Julier2000}. Athans et al. proposed a second-order approximation technique to improve the estimation performance \cite{Athans1968}. This second-order filter requires calculation of both Jacobian and Hessian of the non-linear system under consideration and proved to be computationally expensive. Then N$\phi$rgaard et al. presented an approximate derivative calculation procedure using Stirling's interpolation formula to avoid analytical calculation of the Jacobian and the Hessian for the second-order filter \cite{NøRgaard2000}. This method provides good estimation results for non-linear systems, yet the solution is not exact. A theoretical solution to the non-linear filtering problem requires solving the Fokker-Planck Equation (FPE) which expresses the evolution of the conditional probability density function of the state vector in the form of a partial differential equations \cite{Daum2005}. Daum and Benes discussed the exact solution of non-linear estimation without directly solving the FPE in \cite{Daum1986} and \cite{Benes1985}. However these methods are difficult to implement for high dimensional systems due to computational complexity\cite{Julier2000}. 

 Julier and Uhlmann in their seminal work on non-linear filters, showed a new approach to predict the mean state vector and the error covariance using deterministic sampling \cite{Julier1997,Julier1998,Julier2000,Julier2004}. This approach is known as the Unscented Transformation (UT) and the filter which uses the UT in the prediction step is widely referred as the Unscented Kalman Filter (UKF). The UKF ensures an accuracy of at least the second-order Taylor series approximation without Jacobian and Hessian calculation.  Application of the UKF can be found in satellite navigation, attitude determination and control \cite{Crassidis2003,Fisher2008,Choi2010a}, GPS/INS integration for Unmanned Areal Vehicles \cite{Zhou2010}, indoor positioning \cite{Jiang2014}, target tracking \cite{Leven2009,Li2012} and in various other estimation problems. Unlike the EKF, the UKF does not require an additional process noise matrix and subsequent tuning to compensate for the linearisation. Instead a UKF requires propagation of multiple sampled state vectors which are known as sigma points \cite{Julier2000} to calculate the \textit{a priori} state vector at every time step. For a system with $n$ state elements, $2n+1$ sigma points must be propagated. If the exact difference equation is available for a non-linear system, this propagation of multiple state vectors is not computationally burdensome and the computational effort is comparable to the EKF. But most physical systems are described using non-linear continuous-time differential equations and the system description in the form of difference equations is an approximation of these differential equations. For accurate state propagation, performing a numerical integration is inevitable. From this point of view, multiple numerical integrations must be performed at each UKF prediction stage to calculate the \textit{a priori} state vector, whereas, for the EKF, only one numerical integration operation is required at each step. Due to this reason, implementation of a UKF for a continuous-time system is more computationally expensive than for an EKF. 
 
 The most obvious strategy for reducing the computation time of a UKF is reducing the number of sigma points. From this perspective, several contributions have discussed methods to improve the computational efficiency of the UKF. Julier and Uhlmann showed that at least $n+1$ sigma points are required to capture the uncertainty associated with the system \cite{Julier2002reduced}. The Spherical Simplex Unscented Transform (SSUT) was introduced in their later work. The UKF with SSUT is referred to as the Spherical Simplex Unscented Kalman Filter (SSUKF) in this paper. The SSUT requires $n+2$ sigma points, $n+1$ of which lie on a hypersphere \cite{Julier2003simplex}. However this reduction can lead to a degraded estimation performance \cite{Chang2013MUT} and the reduction in the computational time is intuitively less than 50\% of the same if the UT is used. Chang suggested the Marginal Unscented Transformation (MUT) to reduce the number of sigma points in \cite{Chang2013MUT}. The MUT can be applied to a special type of non-linear function containing linear substructures. It was suggested that, if $n_a$ state elements out of the $n$ state elements are mapped non-linearly then the number of sigma points can be reduced to $2n_a+1$.
 
 Unlike the previously mentioned methods, this paper explores a possibility of reducing the computational time significantly without reducing the number of sigma points. A new approach to the state vector propagation for the UKF prediction stage is presented with the aim of reducing the computational effort. In this method, the \textit{a posteriori} state vector of the previous time step is propagated to the current time step and the other $2n$ sigma points at the current step are calculated approximately from the information of the previous step using a first-order Taylor series approximation. This approach requires propagation of one state vector instead of $2n+1$. The UKF with this strategy is referred as the Single Propagation Unscented Kalman Filter (SPUKF) in the rest of the paper. It is shown in the next section of the paper that the error in the \textit{a priori} state computation using the new method of the order of the second-order terms of the Taylor series expansion and leads to a lower estimation accuracy for the SPUKF compared to the classical UKF. To reduce the error in the \textit{a priori} state prediction, a second method inspired by Richardson Extrapolation is adopted. In the subsequent sections of the paper a UKF with the aforementioned state propagation strategy and the multidimensional extrapolation method is referred to as the Extrapolated Single Propagation Unscented Kalman Filter (ESPUKF). The estimation accuracy and computational effort of the EKF, UKF, SPUKF and the ESPUKF  were compared using a non-trivial non-linear re-entry vehicle tracking problem previously used as a benchmark in \cite{Athans1968,Julier2000,NøRgaard2000} and \cite{Sarkka2007}.
 
 The remainder of the paper is organised as follows: section \ref{sec:prb_stmnt} establishes the problem statement. Section \ref{sec:strat} derives the mathematical formulation for the new state propagation strategy, the error in the \textit{a priori} and the \textit{a posteriori} mean state vectors and error covariances. The extrapolation technique is also derived. In section \ref{sec:complexity} the computation complexity analysis of the SPUKF and the ESPUKF is shown. The two methods are then verified with the re-entry vehicle tracking problem in section \ref{sec:rep_cs}. The state estimation error and computational time for the SPUKF and the ESPUKF are compared with the EKF, UKF and SSUKF. In Section \ref{sec:satnav} the EKF, UKF, SSUKF, SPUKF and ESPUKF are applied to a multi-GNSS based satellite position determination problem and the result is discussed. Section \ref{sec:conclusion} concludes the paper with a discussion of the possibility of real-time application of Unscented Filters using the proposed methods when computational resources are limited.

\section{Problem Statement}
\label{sec:prb_stmnt}
Consider a continuous-time non-linear stochastic dynamical system
\begin{align}
\dot{\bm{Y}}(t) &= {\bm{f}}(t,{\bm{Y}}(t),{\bf \bm{\nu}}(t))\label{eq:system}\\
{\bm{Z}}(k) &= {\bm{h}}({\bm{Y}}(k)) + {\bm{w}}(k)\label{eq:measurement}
\end{align}
%
Here, $t$ denotes continuous time, $k$ is the discrete equivalent of $t$, $\bm{Y}$ is a $n$-dimensional state vector to be estimated from the discrete measurement ${\bm{Z}}(k)$. $\bm{f}$ and $\bm{h}$ are non-linear functions, ${\bf \bm{\nu}}(t)$ is the process noise and ${\bm{w}}(k)$ is the measurement noise. The uncertainty of $\bm{Y}$ is modelled as a probability distribution. In a Kalman Filter framework, the dynamic model of the system is utilized to compute the \textit{a priori} mean and the error covariance of the probability distribution of $\bm{Y}$. To reduce the error in prediction due to model uncertainties, measurement $\bm{Z}$ is used to estimate the \textit{a posteriori} mean state vector and the error covariance of $\bm{Y}$. To apply a classical UKF to this estimation problem, the continuous-time system equation must be converted to a discrete-time difference equation \cite{Julier2000} as shown below.
\begin{equation}
{\bm{Y}}(k+1) = {\bm{F}}(k,{\bm{Y}}(k),{\bf \bm{\nu}}(k))
\label{eq:disc_sys}
\end{equation}
Here $\bm{F}$ is a non-linear function and the discrete time equivalent of $\bm{f}$.  In the UKF framework the UT is used to calculate the predicted mean state vector and the error covariance. In the UT $2n+1$ weighted samples or sigma points are calculated at the $k$th time step from the estimated mean state vector $\widehat{\bm{Y}}^+(k)$ and the error covariance ${\bm{P}}^+_{YY}(k)$ \cite{Julier2000}. The $i$th sigma point $\bm{Y}_i$ and the corresponding weight $W_i$ are \cite{Julier2000}
\begin{align}
{\bm{Y}}_i(k) &= \left \{\begin{array}{ll}
\widehat{\bm{Y}}^+(k) &, (i = 0)\\	
\widehat{\bm{Y}}^+(k) + \Delta {\bm{Y}_i} &, (i = 1,2,3...2n) 
\end{array}\right.\\
W_i &= \left \{\begin{array}{ll}
\frac{\kappa}{n+\kappa}	&, (i = 0)\\
\frac{1}{2(n+\kappa)} &, (i = 1,2,3...2n)
\end{array}\right.
\label{eq:sigmapt_K}
\end{align}
and\\
\begin{tabular}{lcl}
$\Delta{\bm{Y}}_i$ &= $(\sqrt{(n+\kappa){\bm{P}}^+_{YY}(k)})_i$ & for $i = 1,2,3....n$\\
$\Delta{\bm{Y}}_i$ &= $-(\sqrt{(n+\kappa){\bm{P}}^+_{YY}(k)})_i$ & for $i = n+1,2,3....2n$\\
\end{tabular}
Here $(\sqrt{(n+\kappa)\bm{P}^+_{YY}(k)})_i$ is the $i$th column of the matrix $\sqrt{(n+\kappa)\bm{P}^+_{YY}(k)}$.
$\kappa$ is a parameter that can be used for `fine tuning' and if the noise is considered Gaussian, then heuristically $\kappa$ can be selected in such a way that $n+\kappa=3$ \cite{Julier1997}. The sigma points are propagated to the $(k+1)$th time step using (\ref{eq:disc_sys}). The predicted mean $\bm{\widehat{Y}}^-$ at the $(k+1)$th time step can be computed as \cite{Julier2000}
\begin{align}
\bm{\widehat{Y}}^-(k+1) &=\sum_{i=0}^{2n}W_i{\bm{Y}_i}(k+1)\nonumber\\
												&= \sum_{i=0}^{2n}W_i\bm{F}(k,\bm{Y}_i(k),\bm{\nu}(k))
\label{eq:UKF_mean}
\end{align}
and the predicted covariance \cite{Julier2000}
\begin{align}
{\bm{P}^-_{YY}}(k+1)=&\sum_{i=0}^{2n}W_i\left[{\bm{Y}_i}(k+1)-\bm{\widehat{Y}}^-(k+1)\right]\nonumber\\
&\times\left[{\bm{Y}_i}(k+1)-{\bf \widehat{Y}}^-(k+1)\right]^T
\label{eq:UKF_cov}
\end{align}
The predicted measurement is calculated from the computed measurements corresponding to the $i$th propagated sigma point and can be expressed as \cite{Julier2000}
\begin{equation}
{\bm{Z}}_i(k+1) = {\bm{h}}({\bm{Y}}_i(k+1))
\label{eq:sigma_m}
\end{equation}
Then the predicted measurement ${\bm{\widehat Z}}^-$ at $(k+1)$th time step is \cite{Julier2000}
\begin{equation}
{\bm{\widehat Z}}^-(k+1) = \sum_{i=0}^{2n}W_i\bm{Z}_i(k+1)
\label{eq:mean_meas}
\end{equation}
The innovation covariance i.e. the error covariance of the predicted measurement, is \cite{Julier2000}
\begin{align}
\bm{S} = &\sum_{i=0}^{2n}W_i\left[{\bm{Z}}_i(k+1) - {\bm{\widehat Z}}^-(k+1)\right]\nonumber\\
&\times\left[{\bm{Z}}_i(k+1) - {\bm{\widehat Z}}^-(k+1)\right]^T + \bm{R}
\label{eq:innov_cov}
\end{align}
where, $\bm R$ is the measurement noise covariance matrix. The cross covariance matrix is calculated by \cite{Julier2000}
\begin{align}
\bm{P}_{YZ} = &\sum_{i=0}^{2n}W_i\left[{\bm{Y}_i}(k+1)-\bm{\widehat{Y}}^-(k+1)\right]\nonumber\\
&\times\left[{\bm{Z}}_i(k+1) - {\bm{\widehat Z}}^-(k+1)\right]^T
\label{eq:cross_cov}
\end{align}
The measurement update of the estimation can be performed by the standard Kalman Filter equations
\begin{align}
\bm{K} &= {\bm{P}}_{YZ}\bm{S}^{-1}\nonumber\\
\bm{P}^+_{YY} &= \bm{P}^-_{YY} - \bm{K}\bm{S}\bm{K}^T\nonumber\\
\Delta{\bm{Z}} &= {\bm{Z}}(k+1) - {\bm{\widehat Z}}^-(k+1)\nonumber\\
\bm{\widehat{Y}}^+(k+1) &= \bm{\widehat{Y}}^-(k+1) + \bm{K}\Delta{\bm{Z}}\nonumber 
\end{align}
As previously mentioned, the implementation of the UT requires system representation in the form of difference equations. However, the rigorous mathematical models of most physical systems are expressed in the form of differential equations and the difference equation forms of such systems are predominantly first order Taylor series approximations. For example, consider the mathematical process of converting the differential equation (\ref{eq:system}) to a difference equation. (\ref{eq:system}) can be expressed as

\begin{align}
\lim_{\delta t\to 0}\frac{{\bm{Y}}(t+\delta t)-{\bm{Y}}(t)}{\delta t} &= {\bm{f}}(t,{\bm{Y}}(t),{\bf \bm{\nu}}(t))\label{eq:lim} 
\end{align}
This can be approximated as

\begin{align}
{\bm{Y}}(t+\delta t) \approx {\bm{Y}}(t)+{\bm{f}}(t,{\bm{Y}}(t),{\bf \bm{\nu}}(t))\delta t\label{eq:aprx}
\end{align}

If we consider the time interval $\delta t$ in a manner that $t/\delta t = k$ where, $k$ is a positive integer, then the above approximate equation becomes a difference equation

\begin{equation}
{\bm{Y}}(k+1) = {\bm{Y}}(k)+{\bm{f}}(k,{\bm{Y}}(k),{\bf \bm{\nu}}(k))\delta t
\label{eq:approx_eq}
\end{equation}
which is fundamentally a first order Taylor series approximation of ${\bm{Y}}(t)$ around $t$. In the UKF framework, propagation of sigma points requires $2n+1$ evaluations of the function $\bm{f}$. Use of a difference equation for sigma point calculation results in almost the same order of computational efficiency as an EKF depending on the complexity of the system. But this simple approximation can lead to huge prediction error for a highly non-linear system. To avoid this computation error, the differential equation must be used. (\ref{eq:system}) can be rewritten as
\begin{figure}%
\centering
\includegraphics[width=0.95\columnwidth]{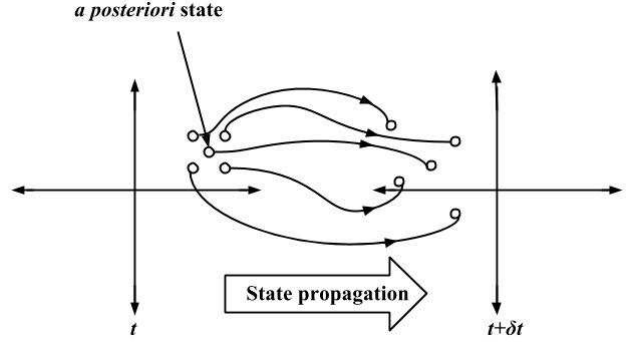}%
\caption{Propagation of sigma points using conventional UT}%
\label{fig:ukf_prop}%
\end{figure}
\begin{align}
{\bm{Y}}(t+\delta t) &= {\bm{Y}}(t)+\int_{t}^{t+\delta t}{{\bm{f}}(\tau,{\bm{Y}}(\tau),{\bf \bm{\nu}}(\tau))}d\tau\nonumber \\
&= {\bm{F}}(t,{\bm{Y}}(t),{\bf\bm{\nu}}(t))\label{eq:propagation}
\end{align}
To determine the \textit{a priori} mean and the error covariance using the UT, $2n+1$ sigma points must be propagated using (\ref{eq:propagation}). The UT approach is shown graphically in Fig. \ref{fig:ukf_prop}. Most of the time a numerical integration technique is adopted for the state propagation using a differential equation, and in this case the computation time for the UT will be much higher. Let us consider a case where a 4th order Runge-Kutta numerical integration method is adopted for the state propagation. Let ${\bm{\widehat{Y}}^+}(t)$ be the \textit{a posteriori} estimate of the state and $\bm{P}^+_{YY}$ be the estimated error covariance at time $t$. Using the Runge-Kutta technique at time $t+\delta t$ the first propagated sigma point will be,
\begin{equation}
{\bm{Y}}_0(t+\delta t) = {\bm{\widehat{Y}}^+}(t) + \frac{\delta t}{6}(\bm{k}_1+2\bm{k}_2+2\bm{k}_3+\bm{k}_4)
\label{eq:runge-kutta}
\end{equation}
where
\begin{align}
\bm{k}_1 &= {\bm{f}}(t,{\bm{\widehat{Y}}^+}(t),{\bf \bm{\nu}}(t))\nonumber\\
\bm{k}_2 &= {\bm{f}}(t+\frac{\delta t}{2},{\bm{\widehat{Y}}^+}(t)+\frac{1}{2}\bm{k}_1\delta t,{\bf \bm{\nu}}(t+\frac{\delta t}{2}))\nonumber\\
\bm{k}_3 &= {\bm{f}}(t+\frac{\delta t}{2},{\bm{\widehat{Y}}^+}(t)+\frac{1}{2}\bm{k}_2\delta t,{\bf \bm{\nu}}(t+\frac{\delta t}{2}))\nonumber\\
\bm{k}_4 &= {\bm{f}}(t+\delta t,{\bm{\widehat{Y}}^+}(t)+\bm{k}_3\delta t,{\bf \bm{\nu}}(t+\delta t))\nonumber
\end{align}
The above-mentioned operations must be performed $2n+1$ times to propagate all the sigma points. As a result, in every prediction step the function $\bm{f}$ must be evaluated $4(2n+1)$ times. However, for accuracy, if $\frac{\delta t}{h}$ $(h \in \mathbb{N})$ step size is chosen to propagate the state vector from time $t$ to $t+\delta t$ then, the function $\bm{f}$ must be evaluated $4h(2n+1)$ times. To account for the process noise in the UKF, the state vector is augmented with the process noise terms \cite{Julier1997},
\begin{equation}
\bm{Y}_a(t) = \begin{bmatrix}
											  \bm{Y}(t)\\
												\bm{\nu}(t)
											\end{bmatrix}	
\label{eq:augmented}
\end{equation}
and the sigma points for the UT are calculated from \cite{Julier1997},
\begin{equation}
\widehat{\bm{Y}}_a(t) = \begin{bmatrix}
													\bm{\widehat Y}^+(t)\\
													\bm{0}_{q\times 1}
													\end{bmatrix}
\label{eq:augmented_0}
\end{equation}
and
\begin{equation}
\bm{P}_a(t) = \left[\begin{array}{ccc}
\bm{P}(t) & \bm{P}_{Y\nu}(t)\\
\bm{P}_{Y\nu}(t) & \bm{Q}(t)
\end{array}\right]
\label{eq:augmented_cov}
\end{equation}
Here $\bm{P}_{Y\nu}$ is the cross covariance matrix of $\bm{Y}(t)$ and $\bm{\nu}(t)$, $q$ is the number of elements in $\bm{\nu}(t)$. The dimension of the augmented state vector is $m = n+q$. Then the number of sigma points to be propagated will be increased to $2m+1$. This puts a substantial computational burden on the processor. Computation time is one of the major constraints of any real-time system. Also, for many real-time applications, the computation power is limited due to several constraints. This makes the UKF difficult to implement in many real-time applications and the EKF remains the first choice of engineers regardless of the fact that UKF can provide a more accurate solution without adding fudge factors in the $\bm{Q}$ matrix. Prior work on improving the computational efficiency of the UKF involves reducing the number of sigma points, compromising the solution accuracy. Therefore, the underlying problem in implementation of the UKF for real-time application is the absence of a sigma point propagation method which is computationally efficient as well as capable of producing an estimation accuracy comparable to the original UKF. In this paper this problem is targeted and a different approach is suggested to significantly improve the computational efficiency of a generic UKF without reducing the number of sigma points.
\begin{figure}%
\centering
\includegraphics[width=0.95\columnwidth]{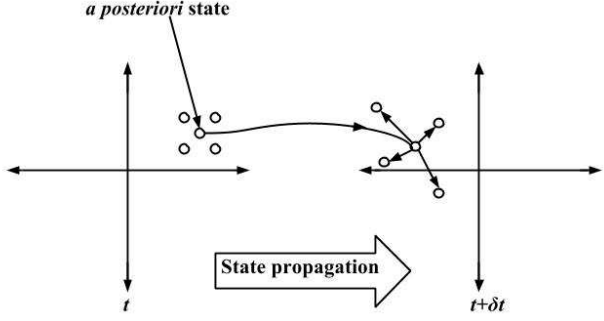}%
\caption{Approximate sigma point propagation}%
\label{fig:ASPP}%
\end{figure}
To reduce the computation time we propose to propagate the \textit{a posteriori} state vector of the previous time step to the current time step using numerical integration as per the EKF and then approximate the other $2n$ (or $2m$) sigma points at the current time step. Fig. \ref{fig:ASPP} shows the hypothesised approach graphically. For approximate calculation of sigma points at the current time step, the deviations of the corresponding sigma points from the \textit{a posteriori} state vector at previous time step  are used in the Taylor series expansion.
\section{Estimation Algorithms Based on the Approximate Sigma Point Propagation}
\label{sec:strat}
\subsection{Single Propagation Unscented Kalman Filter}
\label{sec:proposed_strategy}
The sigma points at time $t$ are
\begin{align*}
\bm{Y}_i(t) &= {\bm{\widehat{Y}}}^+(t)+\Delta\bm{Y}_i, (i = 1,2,3,...2n) 
\end{align*}

The $i$th ($i\neq 0$) sigma point using the UT at time $t+\delta t$ will be

\begin{align}
\bm{Y}_i(t+\delta t) &= \bm{F}(t,\bm{Y}_i(t),{\bm{\nu}}(t))\nonumber\\
& = \bm{F}(t,{\bm{\widehat{Y}}}^+(t),{\bm{\nu}}(t)) + \bm{ D_{\Delta Y_i}F}\nonumber\\
& + \frac{\bm{D_{\Delta Y_i}^2 F}}{2!} + \frac{\bm{D_{\Delta Y_i}^3 F}}{3!}+....
\label{eq:sigmapt}
\end{align}
Here $\bm{ D_{\Delta Y_i}F}=\left.\bm{\frac{\partial F}{\partial Y}}\right|_{\bm{\widehat{Y}}^+(t)}{\Delta \bm{Y}_i}$ is the total differential of $\bm{F}$. If a first order approximation is considered for (\ref{eq:sigmapt}) then 
\begin{align}
\bm{Y}_i(t+\delta t) &\approx \bm{F}(t,{\bm{\widehat{Y}}}^+(t),{\bm{\nu}}(t)) + \bm{ D_{\Delta Y_i}F}\nonumber\\
												&= \bm{Y}_0(t+\delta t) + \left.\bm{\frac{\partial F}{\partial Y}}\right|_{\bm{\widehat{Y}}^+(t)}{\Delta \bm{Y}_i}\label{eq:sigma_approx}
\end{align}
One can easily evaluate $\bm{\frac{\partial F}{\partial Y}}$ as \cite{montenbruck2000satellite}
\begin{equation}
\bm{\frac{\partial F}{\partial Y}} = e^{\bm{\mathcal{J}}\delta t}
\label{eq:}
\end{equation}
where, $\bm{\mathcal{J}}$ is the Jacobian of $\bm{f}$ evaluated at $\bm{Y} = \bm{\widehat{Y}}^+(t)$. By evaluating the Jacobian of $\bm{f}$ once per step, all the sigma points can be propagated over time. Therefore, in this method, instead of evaluating the function $\bm{f}$ $4h(2n+1)$ times, all the sigma points at the current time step are computed by evaluating the function $\bm{f}$ $4h$ times, the Jacobian matrix and matrix exponential once in each step. Section \ref{sec:complexity} will verify that the processing time decreases significantly using this method because of the reduction of the number of evaluation of the function $\bm{f}$. In the SPUKF, this prediction strategy is used. The steps for measurement prediction, Kalman gain computation, the mean state vector and the error covariance calculation is same as the UKF.
\subsubsection{Error in the Mean and the Covariance}
\label{sec:error}
In this subsection, the error incurred due to the approximation mentioned in section \ref{sec:proposed_strategy} is assessed. Hereafter, the sigma points and mean calculated using the proposed approximation will be denoted as $\bm{\widetilde{Y}}_i$ and $\bm{\widetilde{Y}}$ respectively. 

 Let ${\bm{Y}_i}(t)$ be propagated using a 4th order Runge-Kutta numerical integration from (\ref{eq:propagation}). Then, from (\ref{eq:sigmapt}) and (\ref{eq:sigma_approx}) the error due to approximation for each sigma point (neglecting the truncation and round off error) is
\begin{align}
{\bm{Y}_i}(t+\delta t) - {\bm{\widetilde{Y}}_i}(t+\delta t) &= \frac{\bm{D_{\Delta Y_i}^2 F}}{2!} +\frac{\bm{D_{\Delta Y_i}^3 F}}{3!}+....
\label{eq:error_i}
\end{align} 
which is denoted as $\bm{e_i}$. The mean state vector using the original UT \cite{Julier2000}
\begin{align}
{\bm{\widehat{Y}}}^-(t+\delta t) &= \frac{1}{n+\kappa}\left[\kappa{\bm{Y}_0}(t+\delta t) +\frac{1}{2}\sum_{i=1}^{2n}\left[{\bm{Y}_i}(t+\delta t)\right]\right]\nonumber\\
&={\bm{Y}_0}(t+\delta t) + \frac{1}{2(n+\kappa)}\nonumber\\
&\times\sum_{i=1}^{2n}\left[\bm{ D_{\Delta Y_i}F} + \frac{\bm{D_{\Delta Y_i}^2 F}}{2!} +\frac{\bm{D_{\Delta Y_i}^3 F}}{3!}+....\right]\label{eq:mean_UKF}
\end{align}
Similarly the \textit{a priori} mean state vector for the SPUKF is
\begin{align}
\bm{\widetilde{Y}}^-(t+\delta t) &= \frac{1}{n+\kappa}\left[\kappa{\bm{Y}_0}(t+\delta t) +\frac{1}{2}\sum_{i=1}^{2n}\left[{\bm{\widetilde{Y}}_i}(t+\delta t)\right]\right]\nonumber\\
&={\bm{Y}_0}(t+\delta t) + \frac{1}{2(n+\kappa)}\sum_{i=1}^{2n}\bm{ D_{\Delta Y_i}F}\label{eq:approx_ukf}
\end{align}
The error in the prediction due to the approximation is
\begin{align}
&{\bm{\widehat{Y}}}^-(t+\delta t)-\bm{\widetilde{Y}}^-(t+\delta t)\nonumber\\
&= \frac{1}{2(n+\kappa)}\sum_{i=1}^{2n}\left[\frac{\bm{D_{\Delta Y_i}^2 F}}{2!} +\frac{\bm{D_{\Delta Y_i}^3 F}}{3!}+....\right]\label{eq:mean_error}
\end{align}
This mean prediction error is denoted as $\bm{\bar{e}}$.\label{sec:cov}
Let us define the predicted covariance using the original UT as \cite{Julier2000}
\begin{align}
{\bm{P}^-_{YY}}(t+\delta t) &= \frac{1}{n+\kappa}\kappa\left[{\bm{Y}_0}-{\bm{\widehat{Y}}^-}\right]\left[{\bm{Y}_0}-{\bm{\widehat{Y}}^-}\right]^T\nonumber\\
&+\frac{1}{2(n+\kappa)}\sum_{i=1}^{2n}\left[{\bm{Y}_i}-{\bm{\widehat{Y}}^-}\right]\left[{\bm{Y}_i}-{\bm{\widehat{Y}}^-}\right]^T\label{eq:cov}
\end{align}
and the predicted covariance using the new approximation is
\begin{align}
\bm{\widetilde{P}_{YY}}^-(t+\delta t)&=\sum_{i=0}^{2n}W_i\left[{\bm{\widetilde{Y}}_i}-\bm{\widetilde{Y}}^-\right]\left[{\bm{\widetilde{Y}}_i}-{\bm{\widetilde{Y}}}^-\right]^T\nonumber\\
&=\frac{1}{n+\kappa}\kappa\left[{\bm{Y}_0}-{\bm{\widetilde{Y}}^-}\right]\left[{\bm{Y}_0}-{\bm{\widetilde{Y}}^-}\right]^T\nonumber\\
&+\frac{1}{2(n+\kappa)}\sum_{i=1}^{2n}\left[{\bm{\widetilde{Y}}_i}-{\bm{\widetilde{Y}}^-}\right]\left[{\bm{\widetilde{Y}}_i}-{\bm{\widetilde{Y}}^-}\right]^T\label{eq:cov_aprx}
\end{align}
Using (\ref{eq:error_i}), (\ref{eq:mean_error}), (\ref{eq:cov}) and (\ref{eq:cov_aprx}),
\begin{align}
&{\bm{P}^-_{YY}}(t+\delta t)\nonumber\\
&= {\bm{\widetilde{P}}_{YY}}^-(t+\delta t)+\frac{\kappa}{n+\kappa}\left[\bm{\bar{e}}\bm{\bar{e}}^T-\bm{\bar{e}}\left[{\bm{Y}_0}-{\bm{\widetilde{Y}}^-}\right]^T\right.\nonumber\\
&-\left.\left[{\bm{Y}_0} - {\bm{\widetilde{Y}}^-}\right]\bm{\bar{e}}^T\right] + \frac{1}{2(n+\kappa)}\sum_{i=1}^{2n}\left[[\bm{e}_i-\bm{\bar{e}}]\left[{\bm{\widetilde{Y}}_i}-{\bm{\widetilde{Y}}^-}\right]^T\right.\nonumber\\
&+\left.\left[{\bm{\widetilde{Y}}_i}-{\bm{\widetilde{Y}}^-}\right][\bm{e}_i-\bm{\bar{e}}]^T+[\bm{e}_i-\bm{\bar{e}}][\bm{e}_i-\bm{\bar{e}}]^T\right]
\end{align}
For the ease of further calculation, the above equation is rewritten as,
\begin{equation}
{\bm{P}^-_{YY}}(t+\delta t) = {\bm{\widetilde{P}}_{YY}}^-(t+\delta t) + \bm{P}_{YY_e}
\label{eq:Cov_err}
\end{equation}
Here, $\bm{P}_{YY_e}$ is the error in predicted error covariance calculation due to the approximation.
From the above equations it can be inferred that the error in the prediction of the mean and the error covariance with the proposed method is of the order of the second-order Taylor series terms of the original the UT.

\subsubsection{Error in State Estimation}
\label{sec:Error_estimate}
Define $\bm{\widetilde{Z}}_i = \bm{h}(\bm{\widetilde{Y}_i})$ as the computed measurement vector corresponding to the sigma point $\bm{\widetilde{Y}}_i$ and $\bm{\widetilde{Z}}^-$ as the predicted measurement vector at time $t+\delta t$ for the SPUKF. Then,
\begin{align}
\bm{{Z}}_i &= \bm{h}(\bm{Y}_i)\nonumber\\
										 &= \bm{h}(\bm{\widetilde{Y}_i + \bm{e_i}})\nonumber\\
										 &= \bm{\widetilde{Z}}_i + \bm{D}_{\bm{e_i}}\bm{h} + \frac{\bm{D}^2_{\bm{e_i}}\bm{h}}{2!}+....\label{eq:Zi_ZiSP}
\end{align}
Let, $\bm{D}_{\bm{e_i}}\bm{h} + \frac{\bm{D}^2_{\bm{e_i}}\bm{h}}{2!}+.... = \bm{e_{z_i}}$. It is to be noted that, the error term $\bm{e_{z_i}}$is of the order of the first order Taylor series terms, but $\bm{e_i}$ is of the order of the second-order Taylor series terms of (\ref{eq:sigma_approx}). The predicted measurement vector using the UT is
\begin{align}
\bm{\widehat Z}^- &= \sum_{i=0}^{2n}\bm{W}_i\bm{Z}_i\nonumber\\
				 &= \sum_{i=0}^{2n}\bm{W}_i\bm{\widetilde{Z}}_i + \sum_{i=0}^{2n}\bm{W}_i\bm{e_{z_i}}\nonumber\\
				 &= \bm{\widetilde{Z}}^- + \sum_{i=0}^{2n}\bm{W}_i\bm{e_{z_i}}\label{eq:Z_ZSP}
\end{align}
Then, the error in the mean measurement vector prediction using the SPUKF is $\bm{\bar{e}_z} = \sum_{i=0}^{2n}\bm{W}_i\bm{e_{z_i}}$. The innovation $\bm{\Delta Z}$ computed using the UKF can be expressed as
\begin{align}
\bm{\Delta Z}&= \bm{Z} - (\bm{\widetilde Z}+\bm{\bar e}_z)\nonumber\\
						 &= \bm{\Delta\widetilde {Z}} - \bm{\bar e}_z 
\end{align} 
Here, $\bm{\Delta\widetilde {Z}} = \bm{Z} - \bm{\widetilde Z}$ is the innovation computed using the SPUKF.
The cross covariance matrix at time $t+\delta t$ using the UT is
\begin{align}
P_{YZ}&= \sum_{i=0}^{2n}W_i\left[\bm{Y}_i - \bm{\widehat{Y}}^-\right]\left[\bm{Z}_i - \bm{Z}^-\right]^T\nonumber\\
&= \sum_{i=0}^{2n}W_i\left[\bm{\widetilde Y}_i - \bm{\widetilde{Y}}^-\right]\left[\bm{\widetilde Z}_i - \bm{\widetilde Z}^-\right]^T\nonumber\\
&+ \sum_{i=0}^{2n}W_i\left[\bm{e}_i - \bm{\bar{e}}\right]\left[\bm{\widetilde Z}_i - \bm{\widetilde Z}^-\right]^T\nonumber\\
&+ \sum_{i=0}^{2n}W_i\left[\bm{\widetilde Y}_i - \bm{\widetilde Y}^-\right]\left[\bm{e}_{z_i} - \bm{\bar{e}}_z\right]^T\nonumber\\
&+ \sum_{i=0}^{2n}W_i\left[\bm{e}_i - \bm{\bar{e}}\right]\left[\bm{e}_{z_i} - \bm{\bar{e}}_z\right]^T\nonumber\\
&= \bm{\widetilde P}_{YZ} + \bm{P}_{e_{YZ}}\label{eq:P_YZ}
\end{align}
Here, $\sum_{i=0}^{2n}W_i\left[\bm{\widetilde Y}_i - \bm{\widetilde{Y}}^-\right]\left[\bm{\widetilde Z}_i - \bm{\widetilde Z}^-\right]^T = \bm{\widetilde P}_{YZ}$, which is the cross covariance matrix computed using the SPUKF. $\bm{P}_{e_{YZ}}$ is the error in the cross covariance computation. From (\ref{eq:error_i}) and (\ref{eq:mean_error}), it can be deduced that $\bm{P}_{e_{YZ}}$ is of the order of the second-order Taylor series terms. Similarly, the innovation covariance $\bm S$ computed using the UT can also be expressed as,
\begin{equation}
\bm S = \bm{\widetilde S} + \bm{S}_e
\label{eq:S_SSP}
\end{equation}
Here $\bm{\widetilde S}$ is the innovation covariance computed using the SPUKF and $\bm{S}_e$ is the error in innovation covariance calculation, which is also of the order of the second-order Taylor series terms. The Kalman gain $\bm K$ for the original UKF is,
\begin{align}
\bm K &= \bm{P}_{YZ}S^{-1}\nonumber\\
			&= (\bm{\widetilde P}_{YZ} + \bm{P}_{e_{YZ}})(\bm{\widetilde S} + \bm{S}_e)^{-1}\nonumber\\
			&= \bm{\widetilde P}_{YZ}\bm{\widetilde S}^{-1} + \bm{P}_{e_{YZ}}\bm{\widetilde S}^{-1}\nonumber\\
			&- (\bm{\widetilde P}_{YZ} + \bm{P}_{e_{YZ}})[\bm{\widetilde S}^{-1}\bm{S}_e(\bm{S}_e+\bm{S}_e\bm{\widetilde S}^{-1}\bm{S}_e)\bm{S}_e\bm{\widetilde S}^{-1}]\label{eq:K_KSP}
\end{align}
Here, $\bm{\widetilde P}_{YZ}\bm{\widetilde S}^{-1}$ is the Kalman gain $\bm{\widetilde K}$ calculated using the SPUKF. The second and third terms in (\ref{eq:K_KSP}) are the difference in Kalman gains computed using the UKF and the SPUKF and denoted as $\bm{K}_e$. The estimated mean state vector using the UKF is
\begin{align}
\bm{\widehat Y}^+ &= \bm{\widehat Y}^- + \bm{K}\Delta Z\nonumber\\
									&= \bm{\widetilde Y}^- + \bm{\bar e} + (\bm{\widetilde K} + \bm{K}_e)(\bm{\Delta\widetilde{Z}}- \bm{\bar e}_z)\nonumber\\
									&= \bm{\widetilde Y}^- + \bm{\widetilde K}\bm{\Delta \widetilde{Z}} + \bm{\bar e} + \bm{K}_e\bm{\Delta \widetilde{Z}} - \bm{\widetilde K}\bm{\bar e}_z - \bm{K}_e\bm{\bar e}_z\label{eq:estimate_error}
\end{align}
In the above equation $\bm{\widetilde Y}^- + \bm{\widetilde K}\bm{\Delta\widetilde{Z}}$ is the estimated mean state vector using the SPUKF. The fourth, fifth and the sixth terms in the equation can be neglected. This implies that, the error in state estimation in the SPUKF is predominantly due to $\bm{\bar e}$. The estimated error covariance using the UKF is
\begin{align}
\bm{P}^+_{YY} &= \bm{P}^-_{YY} - \bm{K}\bm{S}\bm{K}^T\nonumber\\
							&= \bm{\widetilde P}^-_{YY} + \bm{P}_{YY_e} - [\bm{\widetilde K} + \bm{K}_e][\bm{\widetilde S} + \bm{S}_e][\bm{\widetilde K} + \bm{K}_e]^T\nonumber\\
							&= \bm{\widetilde P}^-_{YY} - \bm{\widetilde K}\bm{\widetilde S}\bm{\widetilde K}^T + \bm{P}_{YY_e}\nonumber\\
							&- [\bm{\widetilde K}\bm{S}_e\bm{\widetilde K}^T + \bm{\widetilde K}\bm{\widetilde S}\bm{K}_e^T + \bm{\widetilde K}\bm{S}_e\bm{K}_e^T + \bm{K}_e\bm{\widetilde S}\bm{\widetilde K}^T\nonumber\\
							&+ \bm{K}_e\bm{S}_e\bm{\widetilde K}^T + \bm{K}_e\bm{\widetilde S}\bm{K}_e^T + \bm{K}_e\bm{S}_e\bm{K}_e^T]\label{eq:estimate_cov_err}
\end{align}
In this equation, $\bm{\widetilde P}^-_{YY} - \bm{\widetilde K}\bm{\widetilde S}\bm{\widetilde K}^T$ is the estimated error covariance using the SPUKF. The remaining terms contributes to the error in the covariance calculation. Except $\bm{P}_{YY_e}$, the other terms are negligibly small. Hence, the error in estimated error covariance using the SPUKF is $\bm{P}_{YY_e}$. From this analysis, it is discerned that, the error in the mean state vector and the error covariance estimation using the SPUKF is of the order of the second-order Taylor series terms compared to the UKF.
\subsection{Multi-dimensional Richardson Extrapolation and the ESPUKF}
\label{sec:richardson1}
As previously discussed, the accuracy of the state estimation using the UKF is  of the order of the second-order Taylor series terms. To achieve similar accuracy using the new sigma point propagation, the second-order Taylor series terms must be included in the \textit{a priori} state vector computation. This can be performed by calculating the Hessian for the function $\bm{f}$. Computation of the Hessian matrix for a complex non-linear system with many state elements is difficult. However, Richardson Extrapolation method is renowned for improving the accuracy of a general approximation technique by an order of the Taylor series terms \cite{richardson1927deferred}. To avoid calculation involving the computation of the Hessian of a non-linear system a time-varying approximation of the Richardson Extrapolation method is adopted to include the second-order Taylor series terms in the \textit{a priori} state vector computation. Using Taylor series approximation,
\begin{align}
{\bm{Y}}_i(t+\delta t) &= {\bm{F}}(t,{\bm{Y}}_i(t),\bm{\nu}(t))\nonumber\\
											  &= {\bm{F}}(t,{\bm{\widehat{Y}}}^+(t),\bm{\nu}(t)) + \bm{ D_{\Delta Y_i}F}|_{{\bm{\widehat{Y}}}^+(t)}\nonumber\\
												&+ \left[\frac{\bm{ D_{\Delta Y_i}^2}}{2!} +\frac{\bm{ D_{\Delta Y_i}^3}}{3!}+....\right]{\bm{F}}|_{{\bm{\widehat{Y}}}^+(t)}\nonumber\\
												&= N_1\bm{ (\Delta Y_i)} + \left[\frac{\bm{ D_{\Delta Y_i}^2}}{2!}+\frac{\bm{ D_{\Delta Y_i}^3}}{3!}+....\right]{\bm{F}}|_{{\bm{\widehat{Y}}}^+(t)}
												\label{eq:method1}
\end{align}
here
\begin{equation}
 N_1\bm{ (\Delta Y_i)} = {\bm{F}}(t,{\bm{\widehat{Y}}}^+(t),\bm{\nu}(t)) + \bm{ D_{\Delta Y_i}F}|_{{\bm{\widehat{Y}}}^+(t)}
\label{eq:N1}
\end{equation}
The Taylor series expansion can also be performed in two steps as shown below:
\begin{align}
&{\bm{Y}}_i(t+\delta t)\nonumber\\ 
&= {\bm{F}}(t,{\bm{\widehat{Y}}}^+(t)+\frac{\bm{\Delta Y}_i}{2}+\frac{\bm{\Delta Y_i}}{2},\bm{\nu}(t))\nonumber\\
											  &= {\bm{F}}(t,{\bm{\widehat{Y}}}^+(t)+\frac{\bm{\Delta Y_i}}{2},\bm{\nu}(t))\nonumber\\
												&+\left[{\bm{D_{\Delta Y_i/2}}}+\frac{\bm{D_{\Delta Y_i/2}^2}}{2!} +\frac{\bm{{D_{\Delta Y_i/2}^3}}}{3!}+....\right]{\bm{F}}|_{{\bm{\widehat{Y}}}^+(t) + \frac{\bm{\Delta Y_i}}{2}}\nonumber\\
												&= {\bm{F}}(t,{\bm{\widehat{Y}}}^+(t),\bm{\nu}(t)) + {\bm{D_{\Delta Y/2}}}{\bm{F}}|_{{\bm{\widehat{Y}}}^+(t)}\nonumber\\
												&+ {\bm{D_{\Delta Y/2}}}{\bm{F}}|_{{\bm{\widehat{Y}}}^+(t) + \frac{\bm{\Delta Y_i}}{2}}\nonumber\\
												&+ \left[\frac{\bm{D_{\Delta Y_i/2}^2}}{2!} +\frac{\bf{{D_{\Delta Y_i/2}^3}}}{3!}+....\right]{\bm{F}}|_{{\bm{\widehat{Y}}}^+(t)}\nonumber\\
												&+ \left[\frac{\bm{D_{\Delta Y_i/2}^2}}{2!} +\frac{\bf{{D_{\Delta Y_i/2}^3}}}{3!}+....\right]{\bm{F}}|_{{\bm{\widehat{Y}}}^+(t) + \frac{\bm{\Delta Y_i}}{2}}\nonumber\\
												&= N_2(\frac{\bm{\Delta Y_i}}{2})+\left[\frac{1}{2^2}\frac{\bm{D^2_{\Delta Y_i}}}{2!} + \frac{1}{2^3}\frac{\bm{D^3_{\Delta Y_i}}}{3!}+...\right]{\bm{F}}|_{{\bm{\widehat{Y}}}^+(t)}\nonumber\\
												&+ \left[\frac{1}{2^2}\frac{\bm{D^2_{\Delta Y_i}}}{2!} + \frac{1}{2^3}\frac{\bm{D^3_{\Delta Y_i}}}{3!}+...\right]{\bm{F}}|_{{\bm{\widehat{Y}}}^+(t)+\frac{\bm{\Delta Y_i}}{2}} \label{eq:method2}
\end{align}
here
\begin{align}
 N_2(\frac{\bm{\Delta Y_i}}{2}) =& {\bm{F}}(t,{\bm{\widehat{Y}}}^+(t),\bm{\nu}(t)) + \bm{D_{\Delta Y/2}}{\bm{F}}|_{{\bm{\widehat{Y}}}^+(t)}\nonumber\\
& + \bm{D_{\Delta Y/2}}{\bm{F}}|_{{\bm{\widehat{Y}}}^+(t) + \frac{\bm{\Delta Y_i}}{2}}
\label{eq:N2}
\end{align}
Considering ${\bm{F}}(t,{\bm{\widehat{Y}}}^+(t),\bm{\nu}(t)) \approx{\bm{F}}(t,{\bm{\widehat{Y}}}^+(t)+\frac{\bm{\Delta Y_i}}{2},\bm{\nu}(t))$ and from (\ref{eq:method1}) and (\ref{eq:method2}) 
\begin{align}
{\bm{Y}}_i(t+\delta t) =& 2N_2(\frac{\bm{\Delta Y_i}}{2})-N_1{\bm{(\Delta Y_i)}}\nonumber\\
&-\left[\frac{1}{2}\frac{\bm{D^3_{\Delta Y_i}}}{3!} + \frac{3}{4}\frac{\bm{D^4_{\Delta Y_i}}}{4!}+...\right]{\bm{F}}|_{{\bm{\widehat{Y}}}^+(t)}
\label{eq:Richardson}
\end{align}
The implication of (\ref{eq:Richardson}) is: if ${\bm{Y}}_i(t+\delta t)$ is approximated using the method $(2N_2(\frac{\bm{\Delta Y_i}}{2})-N_1(\bm{\Delta Y_i}))$, then the error due to the approximation term will be of the order of the 3rd order terms of the Taylor series expansion. With this improvement, a better mean and error covariance propagation is expected and improved estimation accuracy is anticipated. In the ESPUKF, this extrapolation method is utilized to obtain better estimation accuracy. Using the similar approach of section \ref{sec:Error_estimate}, it can be proved that, the error in the state and the error covariance estimation using the ESPUKF is of the order of the third-order Taylor series terms.
\section{Computation Complexity Analysis}
\label{sec:complexity}
To analyse the computation complexity of the UKF, the SPUKF and the ESPUKF, times required to perform mathematical operations are defined as followed:\\ \\
\begin{tabular}{ll}
$t_f$ & time required to evaluate the function $\bm{f}$\\
$t_{mm}$ & time required to multiply a $n\times n$ and\\ 
			& a $n\times n$ matrix\\
$t_m$ & time required to multiply a $n\times n$ and\\ 
			& a $n\times 1$ matrix\\
$t_a$ & time required to add two $n\times 1$ matrix\\
\end{tabular}
\begin{tabular}{ll}
$t_{ms}$ & time required to multiply a scaler with\\ 
				 & a $n\times 1$ matrix\\
$t_{md}$ & time required to divide a $n\times 1$ matrix by a scaler\\
$t_{diag}$ & time required to eigen-decompose a $n\times n$ matrix\\
$t_{sa}$ & time required to perform scaler addition\\
				 & or subtraction\\
$t_{sm}$ & time required to perform scaler multiplication\\
$t_{sd}$ & time required to perform scaler division\\
$t_{exp}$& time required to evaluate scaler exponential\\ 
\end{tabular}\\ \\
Computation complexities of matrix multiplication of two $n\times n$ matrices and matrix multiplication of one $n\times n$ and one $n\times 1$ are $\bm{O}(n^3)$ and $\bm{O}(n^2)$ respectively. Computation complexity of eigen-decomposition is $\bm{O}(n^3)$. Considering $t_{sa}, t_{sm}, t_{sd}, t_{exp}$ as the basic operations, define $t_{max}$ in a fashion that it satisfies the following inequalities,
\begin{align}
t_{max}&\geq max\{t_{sa}, t_{sm}, t_{sd}, t_{exp}\}\label{eq:time_unit}\\
t_{mm} & \leq n^3t_{max}\label{eq:mat_multnn}\\
t_m & \leq n^2t_{max}\label{eq:mat_multn1}\\
t_{diag}& \leq n^3t_{max}\label{eq:diag}
\end{align}
$t_a$, $t_{ms}$ and $t_{md}$ can be expressed as,
\begin{align}
t_a & = nt_{sa}\nonumber\\
		& \leq nt_{max}\label{eq:scaler_add}\\
t_{ms} & = nt_{sm}\nonumber\\
			 & \leq nt_{max}\label{eq:scaler_mult}\\
t_{md} & = t_{sd} + nt_{sm}\nonumber\\
			 & \leq (n+1)t_{max}\label{eq:scaler_div}
\end{align}
If number of basic operations required to evaluate the function $\bm{f}$ is $j (j \in \mathbb{N})$ then,
\begin{equation}
t_f \leq jt_{max}
\label{eq:t_f}
\end{equation}
Here, $j$ reflects the difficulty to evaluate the function $\bm{f}$.
By counting the number of mathematical operations, time required to compute one step of 4th order Runge-Kutta method $t_{RK}$ can be expressed as,
\begin{align}
t_{RK} & = 4t_f + 3t_{sa} + 6t_{ms} + 7t_a + 2t_{sd}\nonumber\\
			 & \leq (13n + 4j + 5)t_{max}\label{eq:t_RK}
\end{align}
Similarly, computation time for Jacobian calculation $t_{\bm{\mathcal{J}}}$ is,
\begin{align}
t_{\bm{\mathcal{J}}} &= nt_f + nt_a + nt_{md} + nt_{sa}\nonumber\\
										 &\leq [2n^2 + (j+2)n]t_{max}\label{eq:t_J}
\end{align}
Computation time to calculate matrix exponential of a $n\times n$ matrix $t_e$ can be expressed as,
\begin{align}
t_e &= t_{diag} + nt_{exp} + 2t_{mm}\nonumber\\
		&\leq (3n^3+n)t_{max}\label{eq:mat_exp}
\end{align}
Computation time required to propagate $2n+1$ state vectors using the UKF is
\begin{equation}
t_{UKF} = h(2n+1)t_{RK}
\label{eq:t_ukf}
\end{equation}
From (\ref{eq:t_RK}),
\begin{align}
t_{UKF}	& \leq [26hn^2 + (8hj+23h)n + 4hj + 5h]t_{max}\label{eq:T_UKF}
\end{align}
here, $h (h \in \mathbb{N})$ is number of steps selected to propagate the state vectors from time $t$ to time $t+\delta t$, i.e. $\delta t/h$ is the step size for the Runge-Kutta algorithm.
 
\subsection{Computation complexity of the SPUKF}
\label{sec:t_SPUKF}
By observing (\ref{eq:sigma_approx}), the computation time required to propagate $2n+1$ state vectors using the SPUKF is
\begin{equation}
t_{SP} = ht_{RK} + t_{\bm{\mathcal{J}}} + t_e + 2nt_m + 2nt_a
\label{eq:T_SP}
\end{equation}
Similar to the UKF, it is assumed that, the Runge-Kutta method uses $\delta t/h$ step size for state propagation.
Using (\ref{eq:mat_multn1}), (\ref{eq:scaler_add}), (\ref{eq:t_RK}), (\ref{eq:t_J}) and (\ref{eq:mat_exp})
\begin{align}
t_{SP} &\leq [5n^3 + 4n^2 + (13h + j + 3)n + 4hj + 5h]t_{max}\label{eq:T_SPUKF}
\end{align}
$t_{UKF}$ and $t_{SP}$ can be written as
\begin{align}
t_{UKF} & = [26hn^2 + (8hj+23h)n + 4hj + 5h]t_{max}\label{eq:T_MAX_UKF}\nonumber\\
				& -\delta_{UKF}~,~\delta_{UKF}\geq 0\\
t_{SP}  & = [5n^3 + 4n^2 + (13h + j + 3)n + 4hj + 5h]t_{max}\label{eq:T_MAX_SP}\nonumber\\
				& -\delta_{SP}~,~\delta_{SP}\geq 0
\end{align}
The terms $\delta_{UKF}$ and $\delta_{SP}$ contain the low level details of the matrix operations if each of the scaler addition, multiplication, division and the scaler exponential operations is considered as one unit operation. Hence, the values of $\delta_{UKF}$ and $\delta_{SP}$ entirely depend on the choice of the algorithms for matrix operations. To make the computation complexity analysis of the Kalman Filters independent of the algorithms used for matrix operations, $\delta_{UKF}$ and $\delta_{SP}$ can be neglected. This assumption leads to comparison of the maximum possible value of $t_{UKF}$ and $t_{SP}$ for a given estimation problem.
 
In the SPUKF framework, only the state propagation method is different to the UKF and the computation time required to generate sigma points at the current time, the computation of the predicted error covariance and the computation of the \textit{a posteriori} mean and error covariance for the SPUKF and the UKF are equal. The reduction in computation time for the SPUKF will be only due to the reduction in computation time for the state propagation. Hence, the difference in computation time between the UKF and the SPUKF is $t_{UKF} - t_{SP}$. For computation time reduction, $t_{UKF} - t_{SP}>0$. The computation time reduction of the SPUKF as a percentage of the computation time of the UKF can be defined as,
\begin{equation}
C_{SP} = \frac{t_{UKF} - t_{SP}}{t_{UKF}}\times 100\%
\label{eq:EFF_SPUKF}
\end{equation}
here, $C_{SP}$ is the percentage computation time reduction for the SPUKF with respect to the UKF. From (\ref{eq:T_MAX_UKF}) and (\ref{eq:T_MAX_SP}) and neglecting $\delta_{UKF}$ and $\delta_{SP}$,
\begin{align}
&C_{SP}\nonumber\\
&= \frac{(26h-4)n^2 + (8hj+10h-j-3)n - 5n^3}{26hn^2 + (8hj+23h)n + 4hj + 5h}\times 100\%\label{eq:C_SP}
\end{align}
\begin{figure}%
\centering
\includegraphics[width=0.95\columnwidth]{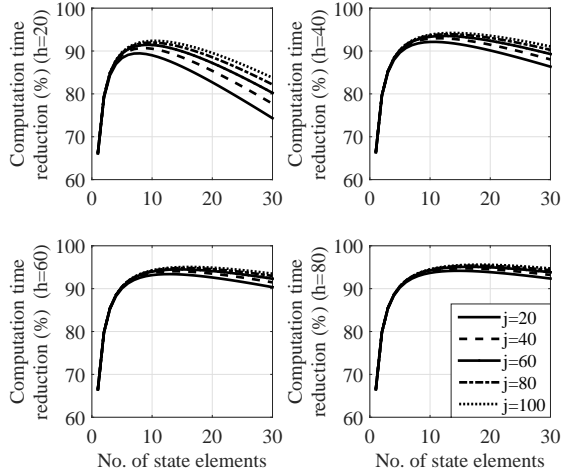}%
\caption{Computation efficiency improvement of the SPUKF}%
\label{fig:eff_SPUKF}%
\end{figure}
\begin{figure}%
\centering
\includegraphics[width=0.95\columnwidth]{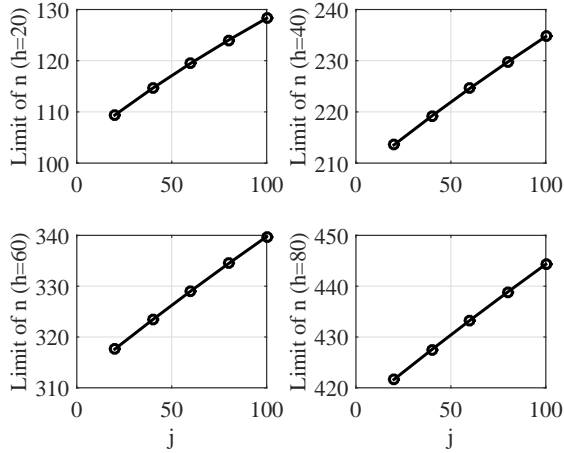}%
\caption{Limit of $n$ for which the SPUKF is more efficient than the UKF}%
\label{fig:lim_SPUKF}%
\end{figure}
This equation implies that, for the computation time of the SPUKF to be reduced compared to the UKF, the numerator of $C_{SP}$ has to be greater than $0$. This limits the number of state elements $n$ that can be estimated using the SPUKF to result in a reduced computation time. Exceeding the limit will result in a higher computation time than the UKF. The limit of $n$ is the real positive root of the polynomial $5n^3 - (26h-4)n^2 - (8hj+10h-j-3)n$,
because at $5n^3 - (26h-4)n^2 - (8hj+10h-j-3)n=0$, $t_{SP} = t_{UKF}$. For different values of $h$ and $j$, $C_{SP}$ vs. $n$ is shown in Fig. \ref{fig:eff_SPUKF}. Fig. \ref{fig:eff_SPUKF} shows that, using the SPUKF more than $90\%$ computation time reduction is possible. With increasing $n$, the computation time reduction increases rapidly and after a certain value of $n$ the reduction decreases slowly. The maximum possible computation time reduction increases with increase in $j$ or $h$. Fig. \ref{fig:lim_SPUKF} shows that, with increasing $j$ and $h$ the limit of $n$ also increases. This implies that, for a higher complexity of the system function $\bm f$ and a higher number of integration steps $h$, the computation time of the SPUKF remains less than that for the UKF at a larger range of $n$, the number of state variables.  

\subsection{Computation complexity of the ESPUKF}
From section \ref{sec:richardson1}, the computation time required to propagate $2n+1$ state vectors using the ESPUKF is
\begin{equation}
t_{ESP} = t_{N1} + 2nt_{N2} + 2nt_{ms} + 2nt_a
\label{eq:T_ESP}
\end{equation}
where, $t_{N1}$ and $t_{N2}$ are the computation time required to calculate $N_1(\Delta Y_i)$  and $N_2(\frac{\Delta Y_i}{2})$. If $\delta t/h$ is the step size chosen for the Runge-Kutta method, then from (\ref{eq:N1}),
\begin{equation}
t_{N1} = t_{SP}\nonumber
\label{eq:T_N1}
\end{equation}
and, from (\ref{eq:N2}),
\begin{align}
t_{N2} & = t_{md} + t_a + t_m + t_{\bm{\mathcal{J}}} + t_e\nonumber\\
			 & \leq [3n^3 + 3n^2 + (j+5)n + j + 1]t_{max}\label{eq:T_N2}
\end{align}
From (\ref{eq:scaler_add}), (\ref{eq:scaler_mult}), (\ref{eq:scaler_div}), (\ref{eq:T_N1}), (\ref{eq:T_N2}) and (\ref{eq:T_ESP})
\begin{align}
t_{ESP} & \leq [6n^4 + 11n^3 + (2j+18)n^2\nonumber\\
				& + (13h+3j+4)n + 4hj +5h]t_{max}
\label{eq:T_ESPUKF}
\end{align}
\begin{figure}%
\centering
\includegraphics[width=0.95\columnwidth]{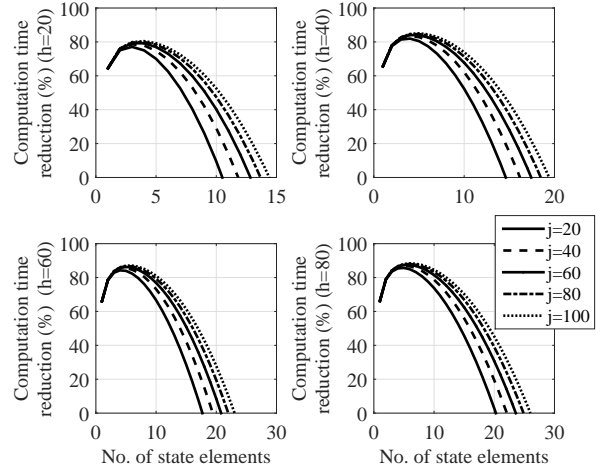}%
\caption{Computation efficiency improvement of the ESPUKF}%
\label{fig:eff_ESPUKF}%
\end{figure}
\begin{figure}%
\centering
\includegraphics[width=0.95\columnwidth]{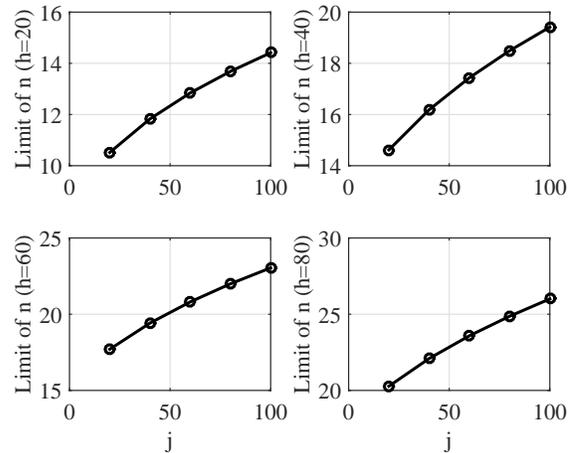}%
\caption{Limit of $n$ for which the ESPUKF is more efficient than the UKF}%
\label{fig:lim_ESPUKF}%
\end{figure}
$t_{ESP}$ can be written as,
\begin{align}
t_{ESP} & = [6n^4 + 11n^3 + (2j+18)n^2\nonumber\\
				& + (13h+3j+4)n + 4hj +5h]t_{max}\nonumber\\
				& -\delta_{ESP}~,~\delta_{ESP}\geq 0\label{eq:T_MAX_ESP}
\end{align}
$\delta_{ESP}$ contains the low level details of the matrix operations and as discussed in \ref{sec:t_SPUKF}, this can be neglected for the purpose of generality. Similar to the SPUKF, the difference in the computation time of the UKF and the ESPUKF is due to the different state propagation strategy. For reduction of the computation time of the ESPUKF, $t_{UKF}-t_{ESP}>0$. The percentage computation time reduction for the ESPUKF with respect to the UKF ($C_{ESP}$) is defined as
\begin{equation}
C_{ESP} = \frac{t_{UKF} - t_{ESP}}{t_{UKF}}\times 100
\label{eq:EFF_ESP}
\end{equation}
From (\ref{eq:T_MAX_UKF}) and (\ref{eq:T_MAX_ESP}) and neglecting $\delta_{UKF}$ and $\delta_{ESP}$,
\begin{align}
&C_{ESP}=\nonumber\\
&\frac{(26h-2k-18)n^2+(8hk+10h-3k-4)n-6n^4-11n^3}{26hn^2 + (8hj+23h)n + 4hj + 5h}\nonumber\\
&\times 100\%\label{C_ESP}
\end{align}
For improvement in computation time, the numerator of $C_{ESP}$ must be greater than $0$. The real positive root of this polynomial $6n^4 + 11n^3 - (26h-2k-18)n^2 - (8hk+10h-3k-4)n$ is the limit of $n$ for the ESPUKF to achieve computation time lower than the UKF. Fig. \ref{fig:eff_ESPUKF} shows $C_{ESP}$ vs $n$ for different values of $j$ and $h$. Depending on the values of $n$, $j$ and $h$, above $80\%$ computation time reduction is possible using the ESPUKF. The characteristics of the graphs are the same as that of the SPUKF however, after reaching the maximum value, the percentage of computation time reduction decreases rapidly with increase in $n$ and as a result the value of $n$ at which, $t_{ESP} = t_{UKF}$ is reduced. Fig. \ref{fig:lim_ESPUKF} shows the limit of $n$ for different values of $j$ and $h$. It is observed that for the ESPUKF at a substantially small value of $n$ than that of the SPUKF, the computation time becomes equal to the computation time of the UKF. This indicates a restricted useful range of the ESPUKF.
\section{Example Application}
\label{sec:rep_cs}
To demonstrate the performance of the SPUKF and the ESPUKF, a non-trivial non-linear re-entry vehicle problem was examined. The problem was used in \cite{Athans1968,Julier2000,NøRgaard2000} and \cite{Sarkka2007} to demonstrate the performance of the respective estimation schemes and as such is a benchmark for this type of work. In the problem a body is considered with a high velocity, which is re-entering the atmosphere at a very high altitude. A radar is used to measure the range of the body in discrete time. The measurement is corrupted by Gaussian noise. The direction of the motion of the body is assumed to be perpendicular to the local horizon. The altitude of the radar is 100000 ft (H) and the horizontal distance between the radar and the body is 100000 ft (M). For convenience of comparison with the work of Julier and Uhlman on the UKF, Imperial units are used in this example. The altitude $x_1(t)$, velocity $x_2(t)$ and the constant ballistic coefficient $x_3(t)$ of the body are to be estimated. The continuous-time dynamics of the system are:

\begin{align}
\dot{x}_1(t) &= -x_2(t) + w_1(t)\\
\dot{x}_2(t) &= -e^{-\lambda x_1(t)}x_2(t)^2x_3(t) + w_2(t)\\
\dot{x}_3(t) &= w_3(t)
\label{eq:dynamics}
\end{align}

here, $w_1(t),w_2(t)$ and $w_3(t)$ are zero-mean, uncorrelated noises with covariance given by ${\bf Q}$ and $\lambda$ is a constant ($5\times 10^{-5}$) that relates the air density and the altitude \cite{Julier2000}. These stochastic dynamical equations are used to generate $x_1$, $x_2$ and $x_3$ and considered as the true value. The random noises used in the generation of the true value of $x_1$, $x_2$ and $x_3$ are considered as unknown during the estimation. All the estimation errors are obtained by subtracting the estimated state element from the corresponding true state element. The range from the radar at time t, $r(t)$ is
\begin{equation}
r(t) = \sqrt{M^2 + [x_1(t)-H]^2} + \nu(t)
\label{eq:measurement_radar}
\end{equation}
here $\nu(t)$ is Gaussian noise. The covariance is considered as $10^4$ ft/s. The initial true state is:
\begin{center}
\begin{tabular}{cl}
$x_1(0)$ &= $30000$ ft\\
$x_2(0)$ &= $20000$ ft/s\\
$x_3(0)$ &= $10^{-3}$
\end{tabular}
\end{center}
The state and the error covariance for the filter initialization are:
\begin{center}
\begin{tabular}{cl}
$\widehat{x}_1(0)$ &= $300000$ ft\\
$\widehat{x}_2(0)$ &= $20000$ ft/s\\
$\widehat{x}_3(0)$ &= $3\times 10^{-5}$\\
\end{tabular}
\end{center}

\begin{align}
\bm{P}(0) &= \left[\begin{array}{ccc}
10^6 & 0 & 0\\
0 & 4\times 10^6 & 0\\
0 & 0 & 10^{-4}
\end{array}\right]\nonumber
\end{align}
The discrete process noise covariance matrix is considered as
\begin{align}
\bm{Q} &= \left[\begin{array}{ccc}
10^{-30} & 0 & 0\\
0 & 10^{-30} & 0\\
0 & 0 & 10^{-30}	
\end{array}\right]\nonumber
\end{align}
\begin{figure}%
\centering
\includegraphics[width=0.95\columnwidth]{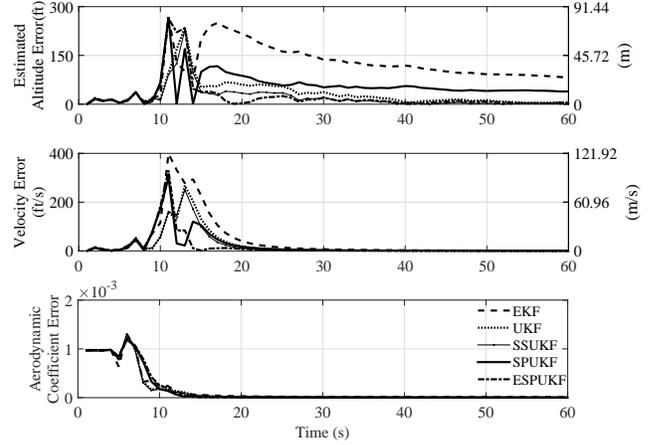}%
\caption{Performance comparison of the SPUKF and ESPUKF with other estimation algorithms}%
\label{fig:EUKF_res}%
\end{figure}

The elements of the $\bm{Q}$ matrix are chosen to be very small value, because arbitrary large values of $\bm{Q}$ matrix entries suppress the approximation error of EKF. For the EKF, UKF, SSUKF, SPUKF and ESPUKF the same initial conditions and $\bm{Q}$ matrix are selected. The state vector is estimated for 1000 seconds using different algorithms separately. Fig. \ref{fig:EUKF_res} shows the estimation errors for the EKF, UKF, SPUKF and ESPUKF for the first 60 seconds of tracking. The altitude errors for the UKF, SSUKF, SPUKF and ESPUKF, are significantly less than for the EKF. But due to the approximate sigma point propagation, the altitude error for the SPUKF is higher than for the original UKF and the SSUKF. The average steady state altitude error and the processing time per time step for different algorithms are presented in Table \ref{tab:EUKF}. The processing time of the SPUKF is significantly lower than for the original UKF and the SSUKF. The average altitude error for the SPUKF is lower than for the EKF. However, the SPUKF provides 14.69 ft (4.48 m) of average altitude error where as the UKF and the SSUKF deliver average altitude errors of 1.19 ft (0.36 m) and 1.53 ft (0.46 m). Estimation using the UKF in general is accurate to the order of the second-order Taylor series terms. From the error analysis presented in sections \ref{sec:error}, \ref{sec:Error_estimate} and the result for the re-entry vehicle tracking problem, the first-order Taylor series approximation for the sigma point propagation is identified as the reason of larger error in the SPUKF.
\begin{figure}[h]%
\centering
\includegraphics[width=0.95\columnwidth]{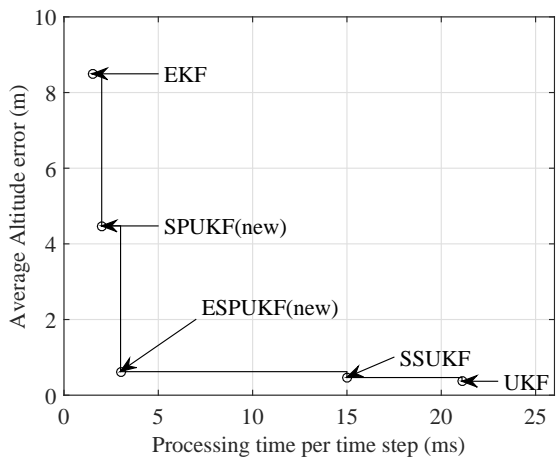}%
\caption{Processing time vs. estimation error for different algorithms}%
\label{fig:comp_graph}%
\end{figure} 
It is observed that the ESPUKF provides lower altitude estimation error compared to the SPUKF due to the elimination of the second-order Taylor series terms. The altitude estimation error using the ESPUKF is comparable with the UKF and the SSUKF results. In Fig. \ref{fig:comp_graph} processing time vs. absolute average steady state estimation error is plotted. It is clear that ESPUKF provides the most desirable performance as compared to the other algorithms when processing time is a constraint.
\section{Application of the SPUKF and the ESPUKF to Satellite Navigation}
\label{sec:satnav}
\begin{table}%
\begin{center}
\caption{Processing time and average altitude error for different algorithms}
\begin{tabular}{|c|c|c|}
\hline
\rule{0pt}{2ex}Algorithm & Processing time & Average steady state\\
& per time step (ms) & altitude error (ft)\\ \hline
\rule{0pt}{2ex}EKF & 1.5 & 27.91\\ 
UKF & 21.1 & 1.19 \\
SSUKF & 14.98 & 1.53 \\
SPUKF & 2.0 & 14.69\\
ESPUKF & 3.0 & 2.04\\ \hline
\end{tabular}
\label{tab:EUKF}
\end{center}
\end{table}
Implementation of a UKF in a real-time on-board satellite navigation scenario was the primary motivation for developing a method of reducing the processing time of the UKF. In the problem under consideration, position of a Low Earth Orbit (LEO) satellite is to be determined from the Global Positioning System (GPS) and Galileo measurements. A SPIRENT Global Navigation Satellite System (GNSS) simulator was used to simulate the the GPS and Galileo constellations, the signals to be received by the user LEO satellite from all the GNSS satellites and the orbital motion of the user LEO satellite. A University of New South Wales (UNSW) Namuru V3.3 multi-GNSS receiver was used to acquire the simulated signal. The pseudo-range measurements are recorded and used to estimate position of the user satellite using the EKF, UKF, SSUKF, SPUKF and the ESPUKF. 
\subsection{System Model}
The state vector associated with the satellite motion is:
\begin{equation}
\bm{X}_{sat} = \begin{bmatrix}
								{\bm{r}}\\
								{\bm{v}}
								\end{bmatrix} = \begin{bmatrix}x &y &z &v_x &v_y &v_z\end{bmatrix}^T
\label{eq:eq1}
\end{equation}
where
${\bm{r}}=[x\quad y\quad z]^T$ is the position vector of the satellite in the Earth Centred Inertial (ECI) frame and
${\bm{v}}=[v_x \quad v_y \quad v_z]^T$ is the velocity vector of the satellite in the ECI frame. $J_2$, $J_3$ and $J_4$ zonal harmonics and neglecting the effect of the Sun and the Moon on the acceleration, the acceleration of the satellite due to the Earth's gravity field are considered in the satellite acceleration model \cite{vallado2001fundamentals,qiao2009multiple}. The differential equation for the states can be represented as:
\begin{equation}
 \bm{\dot{X}}_{sat} = \begin{bmatrix}
												{\bm\dot{r}}\\
												{\bm\dot{v}}
												\end{bmatrix} + \bm{W}_{sat}(t) = \begin{bmatrix}
																														{\bm v}\\
																														{\bm\ddot{r}}
																														\end{bmatrix} + {\bm W_{sat}(t)}
\label{eq:eq3}
\end{equation}
Here the $\bm W_{sat}$ vector is the process noise vector which represents the unmodeled dynamics of the motion. Apart from the states associated with satellite motion, the receiver clock biases for GPS and Galileo are separately considered in the estimable state vector. Although the Namuru V3.3 receiver uses a single clock for both the GPS and Galileo constellations, internally an inter-system bias is introduced. For this reason, two state elements associated with the receiver clock bias are included. A receiver clock bias model can be represented as a first order Markov process. The receiver clock biases dynamics can be represented as \cite{Misra2006}:
\begin{equation}
\bm{\dot{X}}_{clk} = \begin{bmatrix}
											  w_{GPS}\\
												w_{GAL}
											\end{bmatrix}																									
\label{eq:eq4}
\end{equation}
where,
${\bm X_{clk}} = [\delta t_{GPS}\quad \delta t_{GAL}]^T$, ${\bm W_{clk}} = [w_{GPS}\quad w_{GAL}]^T$ is a random noise vector. $\delta t_{GPS}$ is the receiver clock bias for GPS and $\delta t_{GAL}$ is the receiver clock bias for Galileo. The complete state vector and the process noise for the estimation are:
\begin{align*}
\bm{X} = \begin{bmatrix}
						\bm{X}_{sat}\\
						\bm{X}_{clk}
						\end{bmatrix} &&
\bm{W} = \begin{bmatrix}
						\bm{W}_{sat}\\
						\bm{W}_{clk}
						\end{bmatrix} &&
\end{align*}   
\subsection{Measurement model}
Pseudo-range measurements of the GNSS are modelled as \cite{kaplan2005understanding,Misra2006} :
\begin{equation}
\rho_{i}(t) = r_i(t) + c[\delta t_u(t) - \delta t_i(t - \tau)] + I(t) + \epsilon_{\rho}(t)
\label{eq:eq5}
\end{equation}
here, $i$ is the GNSS satellite index, $\rho_i$ is the pseudo-range from the LEO satellite, $\rho_i$ is the pseudo-range from the LEO satellite to the navigation satellite, $r_i$ is the geometric distance from LEO satellite to the navigation satellite, $\delta t_u$ is the receiver clock bias, $\delta t_i$ is the clock bias of the navigation satellite, $\tau$ is the signal transmission time, $c$ is the velocity of light, $I(t)$ is the ionospheric error, $\epsilon_{\rho}(t)$ is the random noise in the pseudo-range measurement. For a space application the ionospheric error is insignificant unless precise position is required. The ionospheric error is neglected in the measurement simulation for this reason. 

\subsection{Simulation Result}
\begin{figure}[h!]%
\centering
\includegraphics[width=0.95\columnwidth]{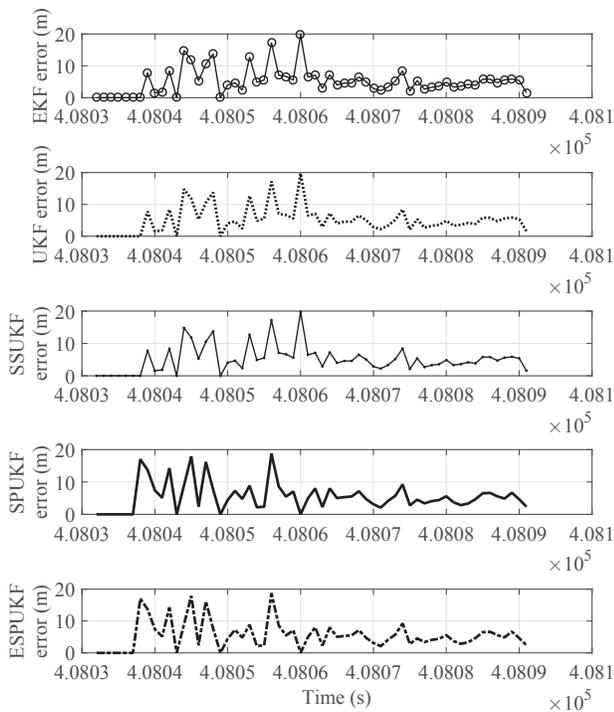}%
\caption{Comparison of the estimated position error}%
\label{fig:sat_filter}%
\end{figure}

Using the system model described by (\ref{eq:eq3}) and (\ref{eq:eq4}), the measurement model equation (\ref{eq:eq5}) and the simulated measurements, the position of the user satellite is estimated using different algorithms. The simulation time span was 120 minutes. The GPS time system was used in the simulation and the time axis in all the plots represents the seconds of the GPS week. The initial position of the satellite was fixed using Least Squares Estimation and used for initialization of all the filters. In Fig. \ref{fig:sat_filter} the norm of the position estimation errors for different estimation algorithms in the ECI frame are plotted for 60 seconds from the starting time. It can be observed that the norm of the estimated position errors lies within $20$m for the EKF, UKF, SSUKF, SPUKF and the ESPUKF. 
\begin{table}%
\caption{Processing Time and Average Position Error for Different Algorithms}
\centering
\begin{tabular}{|l|c|c|}
\hline
\rule{0pt}{2ex}Algorithm & Processing time & Average\\
 &per Time Step (ms) & position error (m)\\ \hline
\rule{0pt}{2ex}EKF & 4.4 & 3.461\\
UKF & 75.9 & 3.448 \\
SSUKF & 45.1 & 3.450 \\
SPUKF & 5.6 & 3.455\\
ESPUKF & 10 & 3.455\\ \hline
\end{tabular}
\label{tab:sat}
\end{table}
Compared to the UKF, the SPUKF and ESPUKF requires much less processing time. The processing time for the SPUKF is reduced by 92.6\% compared to the original UKF. The ESPUKF provides the same average estimation error and the processing time is reduced by 86.8\%. However, the processing time for the ESPUKF is almost twice of the processing time of the SPUKF because the propagated sigma points must be approximated twice. It is to be noted that, the estimation accuracy of the SPUKF and the ESPUKF for the re-entry vehicle problem are different, whereas for the GNSS based satellite navigation problem the estimation accuracies for both the methods are similar. It is anticipated that, for the satellite position estimation problem using the SPUKF the prediction error due to the second-order terms of the Taylor series expansion is not significant but the same has a significant effect on the re-entry vehicle tracking problem. Table \ref{tab:sat} shows the processing time per time step and the estimation errors for different algorithms. 
The computational efficiency of the suggested two methods can be compared with the established processing time reduction techniques like the SSUT and the MUT. State vector prediction using the SSUT requires propagation of $n+2$ sigma points. In this paper it is shown that, theoretically only the \textit{a posteriori} state vector can be propagated to the current time and all the sigma points can be approximated. The results in sections \ref{sec:rep_cs} and \ref{sec:satnav} confirm that the approximation works satisfactorily. In view of the fact that only one state vector is propagated in each time step for the SPUKF and the ESPUKF, it can be argued that the computational efficiency of the aforementioned methods are significantly higher than the SSUKF. The results of the re-entry vehicle tracking problem and the GNSS based satellite navigation problem support the argument. A similar argument can be applied for the comparison with the MUT. Apart from better computational efficiency, the SPUKF and the ESPUKF have another advantage over the MUT. These two methods can be applied to any non-linear system which is continuous in $\mathds{R}^n$ but the MUT can be applied only to the systems which contain linear substructures. 
\section{Conclusion}
\label{sec:conclusion}

	In conclusion, two new methods of implementing unscented estimation algorithms are presented to reduce the computational complexity. Reduction in processing time for the SPUKF and the ESPUKF with respect to the UKF is summarised in table \ref{tab:reduction}.
\begin{table}[h!]%
\caption{Processing time reduction by the new algorithms}
\begin{center}
\begin{tabular}{|c|c|c|}
\hline
\rule{0pt}{2ex}Algorithm & Processing time & Processing time reduction\\
& reduction for re-entry & for LEO satellite\\
&vehicle tracking & position estimation\\ \hline
\rule{0pt}{2ex} SPUKF & 90.5\%  & 92.6\%\\
ESPUKF & 85.5\% & 86.8\%\\ \hline
\end{tabular}	
\end{center}
\label{tab:reduction}
\end{table}
Both the methods have significantly less processing time than the original UKF. This reduction in computation time will make the implementation of the unscented filter easier in a micro-processor with limited computational power for stand-alone real-time application. However, the error characteristics of the two methods are different for the presented two different example applications. These two example applications indicate the existence of two different classes of non-linear estimation problems:
\begin{enumerate}
	\item Non-linear estimation problems for which the UKF provides better accuracy than the EKF. The SPUKF can be used to reduce the computation time and the ESPUKF is required to achieve an estimation accuracy similar to the UKF.
	\item Non-linear estimation problems, for which the UKF and the EKF provide similar accuracy. For this class of problems, the SPUKF and the ESPUKF reduces the computation time for the UKF implementation, however the estimation accuracy is the same for the both algorithms.
\end{enumerate}
From the above observation, it can be hypothesised that there could be a third, intermediate class of non-linear estimation problems, for which the UKF performs better than the EKF and the SPUKF not only reduces the computation time of the UKF but also provides the estimation accuracy similar to the UKF. With further research, a rigorous classification of the non-linear estimation problems can be developed which will facilitate the choice of the estimation algorithm out of the available non-linear estimation techniques.
\bibliographystyle{ieeetr}
\bibliography{biblog}
\section*{Biographies}
\textbf{Sanat Biswas} is a PhD student in the School of Electrical Engineering and Telecommunications at the University of New South Wales (UNSW). He received BE (2010) in Instrumentation and Electronics from Jadavpur University and M. Tech (2012) in Aerospace Engineering from Indian Institute of Technology Bombay. Sanat is currently associated with Australian Centre for Space Engineering Research (ACSER) and Satellite Navigation and Positioning (SNAP) Laboratory. His research focus is non-linear estimation techniques for on-board space vehicle navigation using GNSS receivers. He has been awarded the Emerging Space Leaders Grant 2014 by the International Astronautical Federation.

\textbf{Li Qiao} is a Research Associate in the School of Engineering and Information Technology at the University of New South Wales (UNSW), Canberra. She joined UNSW as a visiting PhD student from 2009 to 2010, and obtained her PhD in Guidance, Navigation and Control at Nanjing University of Aeronautics and Astronautics in 2011. Her research interests are satellite orbit modelling, satellite autonomous navigation and integrated navigation.

\textbf{Andrew Dempster} (M`92\hyp{}SM`03) has a BE (1984) and MEngSc (1992) from UNSW and a PhD from University of Cambridge (1995) in efficient circuits for signal processing arithmetic.

He is Director of the Australian Centre for Space Engineering Research (ACSER) at the University of New South Wales (UNSW). Prof. Dempster was system engineer and project manager for the first GPS receiver developed in Australia in the late 80s and has been involved in satellite navigation ever since. He has published in the areas of arithmetic circuits, signal processing, biomedical image processing, satellite navigation and space systems. His current research interests are in satellite navigation receiver design and signal processing, areas where he has six patents, new location technologies, and space systems.
\end{document}